\documentclass[12pt]{amsart}

\numberwithin{equation}{section}

\newtheorem{theorem}[subsection]{Theorem}
\newtheorem{lemma}[subsection]{Lemma}
\newtheorem{proposition}[subsection]{Proposition}

\theoremstyle{definition}

\newtheorem{definition}[subsection]{Definition}

\theoremstyle{remark}

\newtheorem{remark}[subsection]{Remark}
\newtheorem{remarks}[subsection]{Remarks}
\newtheorem{example}[subsection]{Example}

\newcommand{\F}{{\bf F}}

\newcommand{\SL}{{\text{SL}}}

\newcommand{\LT}{\mathop{\rm LT}}

\newcommand{\LM}{\mathop{\rm LM}}

\newcommand{\binomial}[2]{\genfrac{(}{)}{0pt}{}{#1}{#2}}

\newcommand{\rank}{\mathop{\rm rank}}
\def\Sbar{\overline{S}}
\def\Jbar{\overline{J}}

\title[Resolutions of rings of invariants]{Resolutions of 2 and 3 dimensional rings of invariants for cyclic groups}

\author[Harris]{John C. Harris}
\address{Department of Mathematics and Natural Sciences \\ \hfil\break\indent
				D'Youville College \\ Buffalo, NY, USA \\ 14201
						}
\email{harrisj@dyc.edu}

\author[Wehlau]{David L. Wehlau}
\address{Department of Mathematics and Computer Science \\ \hfil\break\indent
        Royal Military College \\ King\-ston, Ontario, Canada \\ K7K 5L0
             }
\email{wehlau@rmc.ca}

\date{\today}


\dedicatory{}

\begin{document}
\begin{abstract}
 Let $G$ be the cyclic group of order $n$ and suppose ${\bf F}$ is a field containing a primitive $n^\text{th}$ root of unity.
 We consider the ring of invariants ${\bf F}[W]^G$ of a three dimensional representation $W$ of $G$ where $G \subset \text{SL}(W)$.  We describe minimal generators and relations for this ring and prove that the lead terms of the relations are quadratic.   These minimal generators for the relations form a Gr\"obner basis with a surprisingly simple combinatorial structure.  We describe the graded Betti numbers for a minimal free resolution of $F[W]^G$.  The case where $W$ is any two dimensional representation of $G$ is also handled.
\end{abstract}

\maketitle


\section{Introduction}\label{intro}

Let $n \geq 3$ be a positive integer and let $G$ denote the cyclic group of order $n$.
Let $\F$ be a field  which contains a primitive $n^{\rm th}$ root of unity, $\epsilon$.
Let $\sigma$ denote a fixed generator of $G$.
There are $n$ inequivalent irreducible representations of $G$ over $\F$, each of dimension 1.
We denote these by $W_{b}$ with $b=1,2,\dots,n$ where $G$ acts on $W_{b}$ via $\sigma \cdot w = \epsilon^{-b} w$
for all $w \in W_{b}$.

In this paper, we compute the graded Betti numbers for the minimal graded free resolutions for the invariant rings $\F[W]^{G}$, when 
$W=W_{b} \oplus W_{c}$ is two dimensional or $W=W_{b} \oplus W_{c} \oplus W_{d}$ is three dimensional and $G \subset \SL(W)$. 
In each case, we first describe a minimal set of generators for $\F[W]^{G}$ (taken to be monomials since the $G$ action is diagonal) to produce a short exact sequence:
$$0 \to J \to S = \F[x_{1},\dots,x_{m}] \stackrel{\pi}{\longrightarrow} \F[W]^{G} \to 0.$$
By a careful choice of monomial ordering on the ring $S$, we describe a minimal set of homogeneous generators $\{ R_{i,j} \}$ (with $\LT(R_{i,j})=x_ix_j$) for $J$ that is also a Gr\"obner basis for $J$.  In Section~\ref{MonoIdealSec} we find the graded Betti numbers for the monomial ideals $\LT(J)$ in $S$.  It turns out that these minimal graded free resolutions are pure (the differentials have a fixed polynomial degree), so that the graded Betti numbers for $\LT(J)$ and for $J$ are the same.

The two and three dimensional cases are treated in Sections~\ref{2DSec} and \ref{3DSec}, respectively.  
In the three dimensional case, we have $W = W_{b} \oplus W_{c} \oplus W_{d}$, $0 < b,c,d < n$, $G \subset \SL(W)$, and $\F[W] = \F[x,y,z]$.  
Since $G \subset \SL(W)$ we have that $n$ divides $b+c+d$, so that the product of the three linear variables $a = xyz$ is a minimal (since $0 < b,c,d < n$) generator of 
$\F[W]^G$.
It follows that all other minimal generators involve only two of the variables and so correspond to generators arising in the two dimensional case.  
We will describe a set of $m+1$ minimal generators for this ring to produce a short exact sequence:
$$0 \to J \to S = \F[A,B_{1},B_{2},\dots,B_{m}] \stackrel{\pi}{\longrightarrow} \F[W]^{G} \to 0.$$
Also we will describe elements 
$R_{i,j}$  for certain $i,j$ 
with $\LT(R_{i,j})=B_{i}B_{j}$, that form both a minimal generating set for $J$ and also a Gr\"obner basis for $J$.
In the paper \cite{CHW}, an explicit minimal resolution was constructed for this lead term ideal $\LT(J)$ over $S$.  

We would like to thank the referee for a number of expository suggestions and for a simplified proof of Proposition~\ref{IdealBettiProp}.

\section{Minimal Resolutions and Graded Betti Numbers}

 In this section we describe minimal resolutions and some related ideas.  We refer the reader to the books \cite{CLO} and \cite{M-S} for further details.

Let $R = \F[x_1,x_2,\dots,x_m]$ be a graded polynomial ring with each $\deg(x_i)$ a positive integer, and let $I$ be a homogeneous ideal in $R$.  Then the minimal graded free resolution of $R/I$ has the form:
$$0 \rightarrow M_{k} \rightarrow M_{k-1} \rightarrow \cdots \rightarrow M_{1} \rightarrow M_{0}=R \rightarrow R/I \rightarrow 0$$
where $k \leq m$, $M_{i} = \bigoplus_{j} R(-j)^{\beta_{i,j}(R/I)}$, and $R(-j)$ is the free $R$-module obtained by shifting the degrees of $R$ by $j$.  The number $\beta_{i,j}(R/I)$, the $(i,j)$-th {\it graded Betti number} of $R/I$, equals the number of minimal generators of degree $j$ in the $i$-th syzygy module of $R/I$. 

One also defines graded Betti numbers, $\beta_{i,j}(I)$, for the ideal $I$ as follows:
if $\cdots \rightarrow M_{1} \rightarrow M_{0} \rightarrow R/I \rightarrow 0$ is a minimal graded free resolution of $R/I$, then $\cdots \rightarrow M_{1} \rightarrow I \rightarrow 0$ is a minimal graded free resolution of $I$.  Using the standard conventions that $\beta_{-1,0}(I) = 1$ and $\beta_{-1,j}(I) = 0$ for all $j > 0$, it follows that $\beta_{i-1,j}(I) = \beta_{i,j}(R/I)$ for all $i,j \geq 0$.

If we are given a monomial order for $R=\F[x_{1},\dots,x_{m}]$ and if $I$ is a homogeneous ideal in $R$, then we let $\LT(I)$ be the ideal generated by the leading terms of the elements of $I$.  A {\it Gr\"obner basis} for $I$ is a set of elements of $I$ whose lead terms generate $\LT(I)$.  By a simple lead term argument it can be shown that the elements of a Gr\"obner basis for $I$ generate $I$. Since $\LT(I)$ is a monomial ideal, it makes sense to consider the {\it polynomial degree} of elements (defined by setting the degrees of the $x_i$ equal to $1$.)  So for $\LT(I)$ we have two sets of Betti numbers, which we will call the {\it graded} Betti numbers and the {\it polynomial} Betti numbers.

In general, the minimal graded free resolutions of $R/I$ and $R/\LT(I)$ can have different graded Betti numbers.  However, in some cases, they are equal.

\begin{lemma} \label{LeadTermReso} {\rm(\cite{CHW} Section 6)} 
Let $I$ be a homogeneous ideal in the graded ring $R=\F[x_1,\dots,x_{s}]$ minimally generated by a Gr\"obner basis, $\{z_1,\dots,z_{k}\}$, where, for each $i$, the polynomial degree of the lead term of $z_i$ is less than or equal to the polynomial degree of the other terms of $z_i$.  If $(M_*,d_*)$ is a minimal graded free resolution of $R/LT(I)$ that is {\it pure} (i.e., where each of the differentials is homogeneous with respect to the polynomial degree), then there exist differentials $\widehat{d}_*$ so that $(M_*,\widehat{d}_*)$ is a minimal graded free resolution of $R/I$.  In particular, the graded Betti numbers for $R/I$ and $R/\LT(I)$ are the same.
\end{lemma}

In this paper, we are particularly interested in the graded Betti numbers for invariant rings $\F[W]^{G}$.  These are defined as follows: if $\{ b_{1}, \dots, b_{m} \}$ is a minimal set of homogeneous generators for $\F[W]^{G}$, then there is a homomorphism from $R=\F[x_{1},\dots,x_{m}]$ onto $\F[W]^{G}$, where $\deg(x_i) = \deg(b_i)$ and $x_i$ maps to $b_i$.  The kernel of this homomorphism, say $I$, is homogeneous and is called the ideal of relations for $\F[W]^{G}$.  Since $\F[W]^{G} \cong R/I$ as graded rings, we let $\beta_{i,j}(\F[W]^G) = \beta_{i,j}(R/I)$.

It will turn out that Lemma \ref{LeadTermReso} applies to all of the invariant rings considered in this paper.

\section{Two Dimensional Case}\label{2DSec}

  Suppose $W = W_{b} \oplus W_{c}$ is a 2 dimensional representation of $G$ with $0<b,c < n$.
  Then we have $\F[W] = \F[x,y]$ where $\sigma \cdot x = \epsilon^{b} x$ and $\sigma \cdot y = \epsilon^{c} y$.
   It is easy to see that there exist minimal generating sets for the 
   invariant ring $\F[W]^{G}$ consisting of monomials.
   Let $\{u_{i} =  x^{a_i} y^{b_i} \mid i=1,2,\dots,m\}$ be a set of minimal generators:
  $\F[W]^{G} = \F[u_{1},u_{2},\dots,u_{m}]$.  Reordering, if necessary, we may assume that
 $a_{1} > a_{2} > \dots > a_{m-1} > a_{m}$.  Since $m$ is minimal, $u_{i}$ does not divide $u_{j}$  for $i \neq j$.   This implies that
 $b_{1} < b_{2} < \dots < b_{m-1} < b_{m}$.  Moreover, it is clear that $a_1 = n/\gcd(n,b)$, $b_m = n/\gcd(n,c)$ and $a_{m}=0=b_{1}$.

Consider the short exact sequence
$$0 \to J \to S=\F[U_{1},U_{2},\dots,U_{m}] \stackrel{\pi}{\longrightarrow} \F[W]^{G} \to 0$$
where the $U_{i}$ are indeterminants, $\pi(U_{i})=u_{i}$ for $i=1,2,\dots,m$, and $J$ is the ideal of relations for $\F[W]^{G}$.
We put $\deg(U_{i}) := \deg(u_{i})$ and work with the following monomial order on $S$, for monomials $\alpha = \prod_{i=1}^{m} U_{i}^{e_{i}}$ and $\beta = \prod_{i=1}^{m} U_{i}^{f_{i}}$, we say $\alpha < \beta$ if

\begin{enumerate}
  \item $\deg(\alpha) < \deg(\beta)$;
  \item $\deg(\alpha) = \deg(\beta)$ and $e_{i} < f_{i}$ for the largest $i$ where $e_{i} \ne f_{i}$.
\end{enumerate}
For a homogeneous polynomial $f$ in $S$: $\LT(f)$ and $\LM(f)$ will denote the lead term and lead monomial, respectively.

We now describe a minimal generating set for the ideal of relations $J$.
For each pair $i,j$, with $1 \leq i,j \leq m$ and $j-i \ge 2$, consider the product $u_{i}u_{j} \in \F[W]^{G}$.
The monomial $u_{i+1}$ properly divides $u_{i} u_{j}$ and
so $\alpha:=(u_{i}u_{j})/u_{i+1} = x^{a_i+a_j-a_{i+1}} y^{b_i+b_j-b_{i+1}} \ne 1$ is in $\F[W]^{G}$.
Therefore we may write $\alpha=(u_{i}u_{j})/u_{i+1} = \prod_{k=1}^{m} u_{k}^{d_{i,j,k}}$ for some non-negative integers $d_{i,j,k}$.
Note that $u_{r}$ does not divide $\alpha$ for all $r \leq i$ because $a_r \geq a_i > a_i + a_j - a_{i+1}$.  Similarly
$u_{r}$ does not divide $\alpha$ for all $r \geq j$ because $b_r \geq b_j > b_i + b_j - b_{i+1}$.
Therefore $\alpha = \prod_{k=i+1}^{j-1} u_{k}^{d_{i,j,k}}$.
For each pair $i,j$ with $1 \leq i,j \leq m$ and $j-i \ge 2$, fix such a factorization and define
$R_{i,j} := U_{i}U_{j} - U_{i+1}\prod_{k=i+1}^{j-1} U_{k}^{d_{i,j,k}}\in S$.   Note that $R_{i,j} \in J$, that
its lead term $\LT(R_{i,j})=U_{i} U_{j}$ is quadratic, and that its other term has polynomial degree at least $2$.
%
\begin{proposition}\label{Proposition 2D}
The $\binomial{m-1}{2}$ elements
$$\{R_{i,j} \mid 1 \leq i,j \leq m {\rm\ and\ } j-i \geq 2\}$$
form a Gr\"obner basis for $J$ and minimally generate $J$.
\end{proposition}
\proof
We assume, by way of contradiction, that the $R_{i,j}$ do not form a Gr\"obner basis for $J$.
Choose a homogeneous element $f \in J$ with $\beta:=\LM(f)$ minimal such that
$\beta$ is not divisible by $\LT(R_{i,j})$ for all $i,j$.
  We can assume that $f$ is monic, and, since $\F[W]^{G}$
is generated by monomials, that $f$ is a difference of two terms: $f=\prod_{k=1}^{m} U_{k}^{c_{k}}-\prod_{k=1}^{m} U_{k}^{d_{k}}$.
Since $\LT(R_{i,j}) = U_{i}U_{j}$ for $1\leq i,j \leq m$, $j-i \geq 2$ and no leading term divides $\beta$, we conclude that 
$\beta=U_{i}^{p}U_{i+1}^{q}$ for some $1\leq i\leq m$, some positive integer $p$, and some non-negative integer $q$ (where $U_{m+1}:=1$ if $i=m$).
  Recalling that $\beta = \LM(f)$, we have $f = U_{i}^{p}U_{i+1}^{q} - \prod_{k=1}^{i+1} U_{k}^{d_{k}}$ for some non-negative integers $d_k$.
  Since neither $U_{i}$ nor $U_{i+1}$ can divide $f$ (by the minimality of $\beta$)
  we must have  $f  = U_{i}^{p}U_{i+1}^{q} - \prod_{k=1}^{i-1} U_{k}^{d_{k}}$.  We will now show that such an $f$ cannot exist because $U_{i}^{p}U_{i+1}^{q}$ is the smallest monomial in $S$ that maps to $\pi(U_{i}^{p}U_{i+1}^{q}) = u_{i}^{p}u_{i+1}^{q}$.
  Applying $\pi$ to $f$, gives
  $0 = \pi(U_{i}^{p}U_{i+1}^{q}) - \pi( \prod_{k=1}^{i-1} U_{k}^{d_{k}})
         = x^{pa_{i}+q a_{i+1}} y^{pb_{i}+q b_{i+1}} - x^{\sum_{k=1}^{i-1}d_{k}a_{k}}y^{\sum_{k=1}^{i-1}d_{k}b_{k}}$ and therefore
      $pa_{i}+q a_{i+1} = \sum_{k=1}^{i-1}d_{k}a_{k}$ and $pb_{i}+q b_{i+1} = \sum_{k=1}^{i-1}d_{k}b_{k}$.
      Therefore      $$\frac{pb_{i}+q b_{i+1}}{pa_{i}+q a_{i+1}} =  \frac{ \sum_{k=1}^{i-1}d_{k}b_{k}}{\sum_{k=1}^{i-1}d_{k}a_{k}}\ .$$
  Now $p b_i+q b_{i+1} \geq (p+q)b_i$ and $p a_i+q a_{i+1} \leq (p+q)a_i$ and thus
  $$\frac{pb_{i}+q b_{i+1}}{pa_{i}+q a_{i+1}} \geq \frac{(p+q)b_{i}}{(p+q)a_{i}} = \frac{b_{i}}{a_{i}}\ .$$
  Similarly, $\sum_{k=1}^{i-1}d_{k}b_{k} \leq (\sum_{k=1}^{i-1} d_{k}) b_{i-1}$ and
  $\sum_{k=1}^{i-1}d_{k}a_{k} \geq (\sum_{k=1}^{i-1} d_{k}) a_{i-1}$ which implies that
  $$\frac{ \sum_{k=1}^{i-1}d_{k}b_{k}}{\sum_{k=1}^{i-1}d_{k}a_{k}} \leq \frac{(\sum_{k=1}^{i-1} d_{k}) b_{i-1}}{(\sum_{k=1}^{i-1} d_{k}) a_{i-1}} =
     \frac {b_{i-1}}{a_{i-1}}\ .$$
    Since $b_{i}/a_{i} > b_{i-1}/a_{i-1}$ we have a contradiction.   Thus the elements $\LT(R_{i,j})$ give a Gr\"obner basis for $J$ and so also generate $J$.
    
   Finally, since the leading terms of the $R_{i,j}$ are distinct and quadratic, and their non-lead terms have polynomial degree greater than or equal to $2$, they minimally generate. 
     \qed
   \medskip

In Section \ref{MonoIdealSec} below, we will find the polynomial Betti numbers for a class of monomial ideals which includes the ideal $\LT(J)$.
 Proposition~\ref{IdealBettiProp} applied to $\LT(J)$ yields the following result:
\begin{proposition} The polynomial Betti numbers for the ideal $\LT(J)$ are:
$$\beta_{i,j} =
   \begin{cases}
	   1,                                         &{\rm for\ }i=-1,\, j=0;\\
	   (i+1)\binom{m}{i+2} - (m-1)\binom{m-2}{i}, &{\rm for\ }i=0,\dots,m-3,\, j=i+2;\\
	   0,                                         &{\rm otherwise.}
   \end{cases}
   $$
\end{proposition}
Since the non-zero polynomial Betti numbers $\beta_{i,j}$ for a fixed $i$ are concentrated in a single $j$, it follows that the minimal graded free resolution of $\LT(J)$ is pure.  Lemma~\ref{LeadTermReso} then implies that the graded Betti numbers for the ideal $J$ are the same as for $\LT(J)$.  Following the conventions relating Betti numbers of $J$ to Betti numbers of $S/J$, we obtain the following proposition.

\begin{proposition}\label{Minimal Reso 2}
Let $W = W_b \oplus W_c$. Then the minimal free resolution of $\F[W]^G\cong S/J$ as an $S$-module is of the form
$$0 \rightarrow M_{m-2} \stackrel{d_{m-2}}{\longrightarrow} \cdots \rightarrow M_{2} \stackrel{d_{2}}{\rightarrow}
M_{1} \stackrel{d_{1}}{\rightarrow} M_{0} \stackrel{d_{0}}{\rightarrow} \F[W]^G \rightarrow 0,$$
where $\rank(M_{0}) = 1$, and, for $1 \leq i \leq m-2$, $\rank(M_{i}) = i\binom{m}{i+1} - (m-1)\binom{m-2}{i-1}$.  The basis element of $M_{0}$ is in polynomial degree $0$ and the basis elements of $M_{i}$, for $1 \leq i \le m-2$, have polynomial  degree $i+1$.
\end{proposition}

\begin{example}
Let $G={\bf Z}/10$ and let $W = W_b \oplus W_c$ with $b=1$ and $c=2$.  Then $\F[W]^G = \F[x,y]^G$ is generated by the monomials $u_1 = x^{10}$, $u_2 = x^8y$, $u_3 = x^6y^2$, $u_4 = x^4y^3$, $u_5 = x^2y^4$, and $u_6 = y^5$.  The ideal of relations, $J$, has ten minimal generators:
\begin{center}
\begin{tabular}{lll}
$R_{1,3}=U_1U_3-U_2^2$, &$R_{1,4}=U_1U_4-U_2U_3$, &$R_{1,5}=U_1U_5-U_2U_4$,\\ 
$R_{1,6}=U_1U_6-U_2U_5$, &$R_{2,4}=U_2U_4-U_3^2$, &$R_{2,5}=U_2U_5-U_3U_4$,\\
$R_{2,6}=U_2U_6-U_3U_5$, &$R_{3,5}=U_3U_5-U_4^2$, &$R_{3,6}=U_3U_6-U_4U_5$,\\ 
                         $R_{4,6}=U_4U_6-U_5^2$.
\end{tabular}
\end{center}
The polynomial Betti numbers for $S/\LT(J)$ (see Example~\ref{SixExample}) are: \\
\medskip
  \hskip 10pt 
  \begin{tabular}{rrrrrrrr}
      $\beta_{0,0}=1$ &$\beta_{1,2}=10$ &$\beta_{2,3}=20$ &$\beta_{3,4}=15$ &$\beta_{4,5}=4$. \\
  \end{tabular} \\
\medskip
The graded Betti numbers for $S/\LT(J)$ and $S/J$ are: \\
\medskip
	\hskip 10pt
	\begin{tabular}{|r|r|}
			\hline
			           $j:$  &0 \\
			\hline
			$\beta_{0,j}:$ &1 \\
			\hline
	\end{tabular}	
	\hskip 20pt	
	\begin{tabular}{|r|r|r|r|r|r|r|r|}
	    \hline
			           $j:$  &12 &13 &14 &15 &16 &17 &18 \\
			\hline           
			$\beta_{1,j}:$ &1  &1  &2  &2  &2  &1  &1  \\
			\hline
	\end{tabular}\\
	\medskip
	\hskip 10pt	
	\begin{tabular}{|r|r|r|r|r|r|r|r|r|}
	    \hline
			           $j:$  &19 &20 &21 &22 &23 &24 &25 &26 \\
			\hline           
			$\beta_{2,j}:$ &1  &2  &3  &4  &4  &3  &2  &1  \\
			\hline
	\end{tabular}\\ 
	\medskip
	\hskip10pt		
	\begin{tabular}{|r|r|r|r|r|r|r|r|}
	    \hline
			           $j:$  &27 &28 &29 &30 &31 &32 &33 \\
			\hline
			$\beta_{3,j}:$ &1  &2  &3  &3  &3  &2  &1  \\
			\hline
	\end{tabular}
	\hskip 20pt		
	\begin{tabular}{|r|r|r|r|r|}
	    \hline
			           $j:$  &36 &37 &38 &39 \\
			\hline
			$\beta_{4,j}:$ &1  &1  &1  &1  \\
			\hline
	\end{tabular}		
\end{example}

\section{Three Dimensional Case}\label{3DSec}

Now suppose $W = W_{b} \oplus W_{c} \oplus W_{d}$ is a 3 dimensional representation of $G$ with $0< b,c,d < n$ and $G \subset \SL(W)$.
  Then $\F[W] =  \F[x,y,z]$ where $\sigma \cdot x = \epsilon^{b} x$, $\sigma \cdot y = \epsilon^{c} y$, and
$\sigma \cdot z = \epsilon^{d} z$.
Again we can choose a minimal generating set for the $\F$-algebra $\F[W]^{G}$ consisting of monomials.
Since $0< b,c,d < n$ and $G \subset \SL(W)$ (so that $b+c+d \equiv 0 \pmod{n}$),
the element $a = xyz$ is a minimal generator for $\F[W]^G$.
The other generators can be chosen as follows:

$$\bullet\ u_{i} = x^{a_{i}} y^{b_{i}} {\rm\ for\ } i=1,2,\dots,r$$
$${\rm\ with \ } a_{1} > a_{2} > \cdots > a_{r} \ne 0, {\rm\ and\ } 0 = b_{1} < b_{2} < \cdots < b_{r}.$$
$$\bullet\ v_{j} = y^{c_{j}} z^{d_{j}} {\rm\ for\ } j=1,2,\dots,s$$
$${\rm\ with \ } c_{1} > c_{2} > \cdots > c_{s} \ne 0, {\rm\ and\ } 0 = d_{1} < d_{2} < \cdots < d_{s}.$$
$$\bullet\ w_{k} = z^{e_{k}} x^{f_{k}} {\rm\ for\ } k=1,2,\dots,t$$
$${\rm\ with \ } e_{1} > e_{2} > \cdots > e_{t} \ne 0, {\rm\ and\ } 0 = f_{1} < f_{2} < \cdots < f_{t}.$$

\medskip

Consider the short exact sequence
$$0 \to J \to S = \F[A,B_{1},B_{2},\dots,B_{r+s+t}] \stackrel{\pi}{\longrightarrow} \F[W]^{G} \to 0,$$
where $A$ and the $B_{i}$ are indeterminants; $\pi(A)=a$; $\pi(B_{i})=u_{i}$, for $i=1,2,\dots,r$; $\pi(B_{r+j})=v_{j}$, for $j=1,2,\dots,s$; $\pi(B_{r+s+k})=w_{k}$, for $k=1,2,\dots,t$; and $J$ is the ideal of relations.
We set $\deg(A) = \deg(a) = 3$; $\deg(B_{i}) = \deg(\pi(B_{i}))$, for $i = 1,\dots,r+s+t$ and work with the following monomial order on $S$, for monomials
\begin{align*}
\alpha =& A^{g_{0}}  B_1^{g_1} B_2^{g_2} \cdots B_{r+s}^{g_{r+s}}\; B_{r+s+t}^{g_{r+s+1}} B_{r+s+t-1}^{g_{r+s+2}}\cdots B_{r+s+1}^{g_{r+s+t}} \quad\text{and}\\
\beta =& A^{h_{0}}  B_1^{h_1} B_2^{h_2} \cdots B_{r+s}^{h_{r+s}}\; B_{r+s+t}^{h_{r+s+1}} B_{r+s+t-1}^{h_{r+s+2}}\cdots B_{r+s+1}^{h_{r+s+t}},
\end{align*}
we say $\alpha < \beta$ if
\begin{enumerate}
  \item $\deg(\alpha) < \deg(\beta)$,
  \item $\deg(\alpha) = \deg(\beta)$ and $g_{0} > h_{0}$, 
  \item $\deg(\alpha) = \deg(\beta)$, $g_{0} = h_{0}$, and $g_{i} < h_{i}$ for the largest $i$ with $g_{i} \ne h_{i}$. 
\end{enumerate}
(The reader should note the unusual order of the variables in the above expansions of $\alpha$ and $\beta$.)


Note that the three ordered sets $T_{1} = \{B_{1}, B_{2}, \dots, B_{r+1}\}$, $T_{2}=\{B_{r+1}, B_{r+2}, \dots, B_{r+s+1}\}$, and $T_{3}=\{B_{1}, B_{r+s+t}, B_{r+s+t-1}, \dots, B_{r+s+1}\}$ map to the monomials in $\F[W]^{G}$ not divisible by $z$, by $x$, and by $y$, respectively.
This monomial order is carefully chosen to be compatible with these three sets and to ensure that
the elements of the minimal generating set given in Proposition~\ref{Proposition 3D} below form a Gr\"obner basis for $J$.

For the remainder of this section, we will let $m=r+s+t$.

Now we construct a minimal generating set for $J$.   The following distance formula will be useful.  
\begin{definition}
  For $1 \leq i ,j \leq m$ define $|| i-j || := \min\{|i-j|, m-|i-j|\}$.
\end{definition}
This is just the natural distance function on the vertices of a labelled $m$-gon.   

Consider the products $B_{i}B_{j}$ with $1 \leq i<j \leq m$ and $||i - j|| \geq 2$.
If $B_{i}$ and $B_{j}$ are both in $T_{1}$ (they map to monomials not divisible by $z$), or both in $T_{2}$, or both in $T_{3}$,
then construct $R_{i,j}$, with leading term $B_{i}B_{j}$, as in the two dimensional case.
If $B_{i}$ and $B_{j}$ are not both in the same $T_{k}$, then $\pi(B_{i}B_{j})$ is divisible by $a=xyz$, so $\pi(B_{i}B_{j})/a$ is in $\F[W]^G$.  Choose a monomial $\alpha_{i,j} \in S$ with $\pi(\alpha_{i,j})=\pi(B_{i}B_{j})/a$ and define $R_{i,j} = B_{i}B_{j}-A\alpha_{i,j}$, which has leading term $B_{i}B_{j}$ (since a monomial divisible by $A$ is smaller than another of the same degree which is not divisible by $A$).  Of course, for all $i,j$, the non-lead term of $R_{i,j}$ has polynomial degree greater than or equal to $2$.

\begin{proposition}\label{Proposition 3D}
The $\displaystyle\frac{m(m-3)}{2}$ elements
$$\{R_{i,j} \mid  1 \leq i<j \leq m{\rm\ and\ } ||i - j|| \geq 2\}$$
form a Gr\"obner basis for $J$ and minimally generate $J$.
\end{proposition}

\proof
Note that since $\pi$ preserves degree, the ideal $J$ is homogeneous, so can be (minimally) generated by homogeneous elements.  
Also note that if $f$ is any homogeneous element of $S$ for which $A$ divides $\LM(f)$ then $A$ divides $f$.
From this, it follows that $A$ cannot divide the lead monomial of any minimal homogeneous generator of $J$.

We assume, by way of contradiction, that the $R_{i,j}$ do not form a Gr\"obner basis for $J$.
Choose a homogeneous element $f \in J$ with $\beta:=\LM(f)$ minimal such that $\beta$ is not divisible by $\LT(R_{i,j})$ for all $i,j$.
We can assume that $f$ is monic, and, since
$\F[W]^{G}$ is generated by monomials, that $f$ is a difference of two terms:
 $f=\prod_{k=1}^{m} B_{k}^{g_{k}}-A^{h_0}\prod_{k=1}^{m} B_{k}^{h_{k}}$.
Since $\LT(R_{i,j}) = B_{i}B_{j}$ (for $||i-j|| \geq 2$) and no leading term divides $\beta$, we conclude that 
$\beta=B_{i}^{p}B_{i+1}^{q}$ for some $1\leq i\leq m$, some positive integer $p$, and some non-negative integer
$q$ (where $B_{m+1}:=B_{1}$ if $i=m$).  But these possible lead terms involve at most two variables from $\F[W]$, so are the
minimal monomials mapping to $\pi(\beta)$ in $\F[W]^G$ as in the two dimensional case.
Hence these monomials cannot be lead terms for relations and so the $R_{i,j}$ form a Gr\"obner basis for $J$ and so also generate $J$.

Since the leading terms of the $R_{i,j}$ are distinct and quadratic, and their non-lead terms have polynomial degree greater than or equal to $2$, they minimally generate. 
\qed
\medskip

We now describe the minimal graded free resolution of $\F[W]^{G}\cong S/J$.

The following result allows us to work modulo the principle ideal generated by the element $A \in S$.

\begin{lemma} {\rm (\cite{CHW} Lemma 5.3)}
 Let $I$ be a homogeneous ideal in the graded ring $R=\F[x_1,\dots,x_{m}]$ and put $Q=R/I$.
 Suppose $r \in R$ which is not a zero divisor for
$R$ or for $Q$. Let $(M_*,d_*)$ be a minimal graded free resolution of $Q$ as an
$R$-module. Then the complex $(M_*/(r),d_*/(r))$ is a minimal graded free resolution of
$Q/(r)$ as an $R/(r)$-module. In particular, the graded Betti numbers of these
resolutions are the same.
\end{lemma}


The element $A \in S=\F[A,B_{1},\dots,B_{m}]$ is not a zero divisor for $S$ or $S/J$ so the above Lemma applies.
Let $\Sbar = S/(A)$ and $\Jbar = J/(A)$, then the lemma implies that $S/J$ (which is isomorphic to $\F[W]^G$) and $\Sbar/\Jbar$ have the same graded Betti numbers.

In Section \ref{MonoIdealSec}, we will find the polynomial Betti numbers for the lead term ideal $\LT(\Jbar)$, which is the monomial ideal in $\Sbar=\F[B_{1},\dots,B_{m}]$ generated by 
$\{B_{i}B_{j} \mid 1 \leq i < j \leq m {\rm\ and\ } ||i-j|| \geq 2 \}$.  
 Proposition~\ref{IdealBettiProp} applied to $\LT(\Jbar)$ yields the following result:

\begin{proposition}
The polynomial Betti numbers for the ideal $\LT(\Jbar)$ are:
$$\beta_{i,j} =
   \begin{cases}
	   1, &{\rm for\ }i=-1 {\rm\ and\ }j=0;\cr
	   (i+1)\binom{m}{i+2} - m\binom{m-2}{i}, &{\rm for\ }i=0,\dots,m-4 {\rm\ and\ }j=i+2;\\
	   1, &{\rm for\ }i=m-3 {\rm\ and\ }j = m;\\
	   0, &{\rm otherwise.}
   \end{cases}
   $$
\end{proposition}

Since the non-zero polynomial Betti numbers $\beta_{i,j}$ for a fixed $i$ are concentrated in a single $j$, it follows that the minimal graded free resolution of $\LT(\Jbar)$ is pure, so Lemma~\ref{LeadTermReso} implies that the graded Betti numbers for the ideal $\Jbar$ (which are also the same as for $J$) are the same as for $\LT(\Jbar)$.  Again, following the conventions relating Betti numbers of $J$ to Betti numbers of $S/J$, gives the following proposition.
\begin{proposition}\label{Minimal Reso 3}
Let $W = W_b\oplus W_c \oplus W_d$ with $0 < b,c,d < n$ and $b+c+d \equiv 0 \pmod{n}$.
Then the minimal graded free resolution of $\F[W]^G\cong S/J$ as an $S$-module is of the form
$$0 \rightarrow M_{m-2} \stackrel{d_{m-2}}{\longrightarrow} \cdots \rightarrow M_{2} \stackrel{d_{2}}{\rightarrow}
M_{1} \stackrel{d_{1}}{\rightarrow} M_{0} \stackrel{d_{0}}{\rightarrow} \F[W]^G \rightarrow 0,$$
where $\rank(M_{0}) = \rank(M_{m-2}) = 1$, and, for $1 \leq i \leq m-3$, $\rank(M_{i}) = i\binomial{m}{i+1} - m\binomial{m-2}{i-1}$.  The
basis element of $M_{0}$ is in polynomial degree $0$; the basis elements of $M_{i}$, $1 \leq i \le m-3$, have polynomial degree $i+1$; and the basis element of $M_{m-2}$ has polynomial degree $m$.
\end{proposition}

\begin{remark}
A special case of the above result, when $n=1+e+e^2$, $b=1$, $c=e$, and $d=e^2$, can be found in \cite{CHW}. In that paper, we explicitly described basis elements and differentials for the resolution of $\Sbar/\LT(\Jbar)$ and then constructed a contracting homotopy to prove exactness.  The techniques used in Section~\ref{MonoIdealSec} are much simpler.
\end{remark}

\begin{example}
Let $G={\bf Z}/6$ and let $W = W_b \oplus W_c \oplus W_d$ with $b=1$, $c=2$, and $d=3$.  Then $\F[W]^G = \F[x,y,z]^G$ is generated by the monomials $a=xyz$, $b_1=x^6$, $b_2=x^4y$, $b_3=x^2y^2$, $b_4=y^3$, $b_5=z^2$, and $b_6=x^3z$.  The ideal of relations, $J$, has nine minimal generators: 
\begin{center}
\begin{tabular}{lll}
  $R_{1,3}=B_1B_3-B_2^2$, &$R_{1,4}=B_1B_4-B_2B_3$, &$R_{1,5}=B_1B_5-B_6^2$,\\
  $R_{2,4}=B_2B_4-B_3^2$, &$R_{2,5}=B_2B_5-AB_6$,   &$R_{2,6}=B_2B_6-AB_1$,\\
  $R_{3,5}=B_3B_5-A^2$,   &$R_{3,6}=B_3B_6-AB_2$,   &$R_{4,6}=B_4B_6-AB_3$.\\
\end{tabular}
\end{center}
\medskip
The polynomial Betti numbers for $S/\LT(J)$ (see Example~\ref{SixExample}) are: \\
\medskip
  \hskip 10pt 
  \begin{tabular}{rrrrrrrr}
      $\beta_{0,0}=1$ &$\beta_{1,2}=9$ &$\beta_{2,3}=16$ &$\beta_{3,4}=9$ &$\beta_{4,6}=1$. \\
  \end{tabular} \\
\medskip
The graded Betti numbers for $S/\LT(J)$ and $S/J$ are: \\
\medskip
{\footnotesize
	\begin{tabular}{|r||r|r|r|r|r|r|r|r|r|r|r|r|r|r|r|r|r|}
	    \hline
			          $j:$ &0 &$\dots$ &6 &7 &8 &9 &10 &11 &12 &13 &14 &15 &16 &17 &18 &$\dots$ &24 \\
			\hline
			\hline
			$\beta_{0,j}:$ &1 &        &  &  &  &  &   &   &   &   &   &   &   &   &   &        &   \\
			\hline
			$\beta_{1,j}:$ &  &        &1 &2 &3 &2 &1  &   &   &   &   &   &   &   &   &        &   \\
			\hline
			$\beta_{2,j}:$ &  &        &  &  &  &  &2  &4  &4  &4  &2  &   &   &   &   &        &   \\
			\hline
			$\beta_{3,j}:$ &  &        &  &  &  &  &   &   &   &   &1  &2  &3  &2  &1  &        &   \\
			\hline
			$\beta_{4,j}:$ &  &        &  &  &  &  &   &   &   &   &   &   &   &   &   &        &1  \\
			\hline
	\end{tabular}	
}
\end{example}

\section{Graded Betti Numbers for Some Monomial Ideals}\label{MonoIdealSec}

Let $R=\F[x_{1},\dots,x_{m}]$, and let $I$ be an ideal in $R$ that is generated by square-free quadratic monomials.  Then $I$ can be described as the {\it edge ideal} of a simple graph.  Here are the relevant definitions and results we will need (we follow the exposition of \cite{HvT}).

Let $X$ be a finite graph with vertex set $V_{X} = \{x_1,x_2,\dots,x_n\}$ and edge set $E_{X}$. The graph $X$ is called {\it simple} if $X$ has no loops or multiple edges.  A simple graph need not be connected - we let $\#{\rm comp}(X)$ denote the number of components of $X$.  When $Y$ is a subset of $V_{X}$, the {\it induced subgraph} of $X$ on the vertex set $Y$, denoted $X_{Y}$, is the subgraph of $X$ with vertices $Y$ and edge set consisting of the edges in $X$ connecting vertices in $Y$.  The {\it complement} of a graph $X$, denoted $X^{c}$, is the graph whose vertex set is the same as $X$, but whose edge set is the complement of the edge set of $X$: the edge $\{x_{i},x_{j}\} \in E_{X^{c}}$ if and only if $\{x_{i},x_{j}\} \not\in E_{X}$. The {\it complete graph} on a set of vertices is the graph with exactly one edge for each pair of distinct vertices.

A {\it simplicial complex} $\Delta$ on a vertex set $V_{\Delta}$ is a collection of subsets of $V_{\Delta}$ such that: (i) for each vertex $x_{i} \in V_{\Delta}$, the set $\{x_{i}\} \in \Delta$, and (ii) for each set $F$ in $\Delta$, if $E \subseteq F$, then $E \in \Delta$.
An element $F$ of a simplicial complex $\Delta$ is also called a {\it face} of $\Delta$.
The dimension of a face $F$ of $\Delta$, denoted $\dim{F}$, is defined to be $|F|-1$, where $|F|$ denotes the number of vertices in $F$.
The maximal faces of $\Delta$ under inclusion are called the {\it facets} of $\Delta$.
There is a one-to-one correspondence between simplicial complexes and their facet sets.

If $\Delta$ is a simplicial complex with vertex set $V_{\Delta} = \{ x_{1},\dots,x_{m} \}$, then we can associate
two ideals to $\Delta$ in the polynomial ring $R = \F[x_{1},\dots,x_{m}]$. (By abuse of notation we use $x_{i}$ to denote both a vertex of $\Delta$ and a variable in $R$.)  For a face $F$ of $\Delta$, we let $x^{F}$ be the monomial $\prod_{x \in F} x$ in $R$.  The {\it facet ideal}, $\mathcal{I}(\Delta)$, and {\it Stanley-Reisner} ideal, $I_{\Delta}$, of $\Delta$ are defined as follows:
$$\mathcal{I}(\Delta) = \langle x^{F} \mid F {\rm \ is\ a \ facet\ of\ }\Delta \rangle {\rm\ and\ }
I_{\Delta} = \langle x^{F} \mid F \subseteq V_{\Delta} {\rm \ and\ } F \not\in \Delta \rangle.$$

Note that a connected graph $X$ can be thought of as a simplicial complex with facet set equaling its edge set. In this situation, the facet ideal $\mathcal{I}(X)$ is also called the {\it edge ideal}.  Additionally, we can associate another simplicial complex to $X$: the {\it clique complex} of $X$ is the simplicial complex $\Delta(X)$ where $F = \{x_{i_{1}},\dots,x_{i_{j}}\} \in \Delta(X)$ if and only if the induced subgraph $X_{F}$ is a complete graph.  Note that when $X$ is a simple graph, then the edge ideal $\mathcal{I}(X)$ is generated by square-free quadratic monomials and is equal to the Stanley-Reisner ideal of the clique complex of its complement, that is, $\mathcal{I}(X) = I_{\Delta(X^{c})}$.  The following result uses Hochster's and Eagon-Reiner's formulas to compute the graded Betti numbers of $\mathcal{I}(X)$.

\begin{theorem}\label{BettiX} {\rm (\cite{HvT}, Theorem 3.2.1)}
Let $X$ be a simple graph with edge ideal $\mathcal{I}(X)$.  Then
$$\beta_{i,j}(\mathcal{I}(X)) = \sum_{Y \subseteq V_{X},\, |Y|=j} \dim_{\F} \widetilde{H}_{j-i-2} (\Delta(X^{c}_{\, Y}),\, \F)\ for\ all\ i,j \geq 0.$$
\end{theorem}

In the case when $j = i+2$, the above formula involves only the $0$-th homology, so we also have:

\begin{theorem}\label{BettiIplusTwo} {\rm (\cite{HvT}, Theorem 3.2.4)}
Let $X$ be a simple graph with edge ideal $\mathcal{I}(X)$.  Then for all $i \ge 0$,
$$\beta_{i,i+2}(\mathcal{I}(X)) = \sum_{Y \subseteq V_{X},\, |Y|=i+2} (\#{\rm comp}(X^{c}_{\,Y}) - 1).$$
\end{theorem}

We now apply these results to compute the Betti numbers for $\mathcal{I}(X)$ for a specific collection of graphs $X$. 

Let $X[0]$ be the complete graph with vertices $V_{X} = \{ x_{1},\dots,x_{m} \}$ and all $\binom{m}{2}$ edges. The $m$ edges $\{\, \{x_{i},x_{i+1}\} \mid 1 \le i \le m \}$ (where $x_{m+1}$ will mean $x_{1}$) will be called {\it adjacent edges}.
Let $X[s]$ be some graph obtained from the complete graph $X[0]$ by deleting $s$ of the adjacent edges where $1 \leq s \leq m$.

\begin{proposition}\label{IdealBettiProp}
For $s=0,\dots,m-1$, the graded Betti numbers of $\mathcal{I}(X[s])$ are:
$$\beta_{i,j} =
   \begin{cases}
	   1, &{\rm for\ }i=-1 {\rm\ and\ }j=0;\cr
	   (i+1)\binom{m}{i+2} - s\binom{m-2}{i}, &{\rm for\ }i=0,\dots,m-2 {\rm\ and\ }j=i+2;\\
	   0, &{\rm otherwise.}
   \end{cases}
   $$
For $s=m$, the graded Betti numbers of $\mathcal{I}(X[m])$ are:
$$\beta_{i,j} =
   \begin{cases}
	   1, &{\rm for\ }i=-1 {\rm\ and\ }j=0;\cr
	   (i+1)\binom{m}{i+2} - m\binom{m-2}{i}, &{\rm for\ }i=0,\dots,m-4 {\rm\ and\ }j=i+2;\\
	   1, &{\rm for\ }i=m-3 {\rm\ and\ }j = m;\\
	   0, &{\rm otherwise.}
   \end{cases}
   $$

\end{proposition}

\proof
  By definition, $\beta_{-1,0}=1$ and $\beta_{-1,j}=0$ if $j \neq 0$.
  
  For $i \geq 0$, we use Theorem~\ref{BettiX} and Theorem~\ref{BettiIplusTwo}.  It is easy to see that $X[m]^c$ is just an $m$-cycle.  It is also easy to see that $X[s]^c$ is a subgraph of this $m$-cycle with the same vertices but only $s$ edges.    
  Hence for any subset $Y \subseteq V_{X}$, each component of the clique complex $\Delta(X[s]^{c}_{\, Y})$ is simply connected unless $s=m=|Y|$. 

  First consider the case where $j=|Y| < m$ so that 
  $$\widetilde{H}_{j-i-2} (\Delta(X[s]^{c}_{\, Y}),\, \F)=0$$
   unless $j-i-2=0$, and consequently, $\beta_{i,j} = 0$ unless $j=i+2$.
   We use Theorem~\ref{BettiIplusTwo} to compute the remaining Betti numbers for $j=i+2$.     Since any acyclic graph with $j$ vertices and $t$ edges has $j-t$ components, we
   see that 
   \begin{align*}
   \beta_{j-2,j}(\mathcal{I}(X[s])) &= \sum_{Y \subseteq V_{X},\, |Y|=j} (\#{\rm comp}(X[s]^{c}_{\,Y}) - 1)\\
                                                  &= \sum_{Y \subseteq V_{X},\, |Y|=j} j-1-\#E_{X[s]^{c}_{\, Y}}\\
                                                  &= (j-1)\binom{m}{j} - \sum_{Y \subseteq V_{X},\, |Y|=j} \#E_{X[s]^{c}_{\, Y}}\ .
   \end{align*}
     Since each of the $s$ edges of $X[s]^{c}$ lies in exactly $\binom{m-2}{j-2}$ subsets $Y\subseteq V_{X}$ of size $j$
   we see that 
   $$
   \beta_{j-2,j}(\mathcal{I}(X[s])) =  (j-1)\binom{m}{j}  - s\binom{m-2}{j-2} = (i+1)\binom{m}{i+2} - s\binom{m-2}{i}
   $$
   for $i=j-2$ as required.



  Finally, we consider the case where $j=|Y|=m$ (so necessarily, $s=m$ and $Y=V_X$) and $i \geq 0$.   By Theorem~\ref{BettiX} we have $\beta_{i,m}(\mathcal{I}(X[m])) = \dim_{\F} \widetilde{H}_{m-i-2} (\Gamma,\, \F)$, where 
  $\Gamma = \Delta(X[m]^{c}_{\, Y}) = \Delta(X[m]^c)$ is the clique complex of the $m$-cycle, i.e., $\Gamma$ is a circle.  Thus $\beta_{m-3,m}(\mathcal{I}(X[m]))=1$ and $\beta_{i,m}(\mathcal{I}(X[m]))=0$ 
  for $i \neq m-3$.  
\qed
\medskip

Note that the lead term ideals in Sections \ref{2DSec} and \ref{3DSec} correspond to the edge ideals $\mathcal{I}(X[m-1])$ and $\mathcal{I}(X[m])$, respectively.  Also note that the resolutions for $\mathcal{I}(X[s])$, $s=0,\dots,m-1$ are {\it linear} (see \cite{RvT} for the definition of a linear resolution).

\begin{example}\label{SixExample} Here are the polynomial Betti numbers, $\beta_{i,j}$, for the $\mathcal{I}(X[s])$ when $m=6$.  The first column corresponds to $j=0$, the last non-zero entry (for $\beta_{3}(\mathcal{I}(X[6])$) corresponds to $j=6$, and all other entries correspond to $j=i+2$.
\end{example}
\hskip 30pt
		\begin{tabular}{|c||r|r|r|r|r|r|r|}
			\hline
			$i=$ &\ -1 &\ \ 0 &\ \ 1 &\ \ 2 &\ \ 3 &\ \ 4 \\
			\hline
			\hline
			$\beta_{i}(\mathcal{I}(X[0]))$ &1 &15 &40 &45 &24 &5 \\
			\hline
			$\beta_{i}(\mathcal{I}(X[1]))$ &1 &14 &36 &39 &20 &4 \\
			\hline
			$\beta_{i}(\mathcal{I}(X[2]))$ &1 &13 &32 &33 &16 &3 \\
			\hline
			$\beta_{i}(\mathcal{I}(X[3]))$ &1 &12 &28 &27 &4  &2 \\
			\hline
			$\beta_{i}(\mathcal{I}(X[4]))$ &1 &11 &24 &21 &8  &1 \\			
			\hline
			$\beta_{i}(\mathcal{I}(X[5]))$ &1 &10 &20 &15 &4  &0 \\
			\hline
			$\beta_{i}(\mathcal{I}(X[6]))$ &1 &9  &16 &9  &1  &0 \\
			\hline
	\end{tabular}		
\medskip


\begin{remark} One might notice that the final row of this table is symmetric.  In the application to rings of invariants in Section \ref{3DSec}, this corresponds to the fact that the invariant ring $\F[W]^G$ is Gorenstein when $G \subseteq {\rm SL}(W)$ (see \cite{Wat}) so its minimal graded free resolution is self-dual. 
\end{remark}

\begin{remarks}
Minimal resolutions for $\mathcal{I}(X[0])$ and $\mathcal{I}(X[m])$ were described in \cite{CHW}, where explicit basis elements and differentials were given and the complex was proved to be exact using a contracting homotopy.  Also, the Betti numbers for $\mathcal{I}(X[0])$ were found in \cite{J}, Section 5, and in \cite{RvT}, Example 2.10.  The Betti numbers for $\mathcal{I}(X[s])$, for $s=0$, $m-1$, and $m$, go back even further to a paper on rational singularities by Wahl \cite{W}.
\end{remarks}


\ifx\undefined\bysame
\newcommand{\bysame}{\leavevmode\hbox to3em{\hrulefill}\,}
\fi


\begin{thebibliography}{CHW}

%
%
%
\bibitem[CHW]{CHW}
H.E.A.~Campbell, J.C.~Harris and David L.~Wehlau,
{\em Internal Duality for Resolutions of Rings},
J.~of Algebra  {\bf 215} (1999), 1--33.
%
%
\bibitem[CLO]{CLO}
D.~Cox, J.~Little, and D.~O'Shea,
 {\em Ideals, varieties, and algorithms},
Springer-Verlag, 1992.
%
%
%
\bibitem[HvT]{HvT}
H.~T.~Ha and A.~Van Tuyl,
 {\em Resolutions of Square-Free Monomial Ideals Via Facet Ideals: A Survey},
arXiv:math/0604301v2 [Math.AC] 28 May 2007.
%
\bibitem[J]{J}
S.~Jacque,
 {\em Betti Numbers of Graph Ideals},
arXiv:math.AC/0410107 v1, 5 October 2004.
%
\bibitem[M-S]{M-S}
E. Miller and B. Sturmfels, Combinatorial Commutative Algebra. Springer GTM {\bf 227}, Springer, 2004.
%
\bibitem[RvT]{RvT}
M.~Roth and A.~Van Tuyl,
 {\em On the Linear Strand of an Edge Ideal},
arXiv:math.AC/0411181 v2, 12 April 2006.
%
\bibitem[Wah]{W}
J.~M.~Wahl,
{\em Equations defining rational singularities},
Ann.~Sci.~\'Ecole Norm.~Sup. {\bf 10} (1977), 231--264.
%
\bibitem[Wat]{Wat}
K.~Watanabe,
{\em Certain invariant subrings are Gorenstein. I},
Osaka J.~Math. {\bf 11}, No.~1 (1974), 1--8. 
\end{thebibliography}
\end{document}